\documentclass[twocolumn]{revtex4}
\pdfoutput=1
\usepackage{latexsym}
\usepackage{amsmath}
\usepackage{times}
\usepackage{amssymb}
\usepackage{fancyheadings}
\usepackage[T1]{fontenc}
\usepackage[utf8]{inputenc}
\usepackage{url}
\usepackage{hyperref}
\usepackage{color}
\usepackage{graphicx}
 \usepackage{epsfig}
\usepackage{mathptmx}
\usepackage{bm}
\usepackage{array}
\usepackage{algorithm}
\usepackage{algorithmic}
\usepackage{url}
\usepackage{hyperref}
\usepackage{algorithm}
\usepackage{algorithmic}
\usepackage{color}

\usepackage{lscape}
\usepackage{draftcopy}

\begin{document}

\title{\bf Topology of nanonetworks grown by aggregation of simplexes
  with defects}
\author{Bosiljka
  Tadi\'c$^{1,6}$, Milovan \v{S}uvakov$^{2,5}$, Miroslav
  Andjelkovi\'c$^{3}$, Geoff J. Rodgers$^{4}$}
\vspace*{3mm}
\affiliation{$^1$Department of Theoretical Physics, Jo\v zef Stefan Institute,
Jamova 39, Ljubljana, Slovenia}
\affiliation{$^2$Institute of Physics,  University of Belgrade, Pregrevica 118, 11080
  Zemun-Belgrade, Serbia}
\affiliation{$^3$Institute for Nuclear Sciences Vin\v{c}a, University of Belgrade, 11000 Belgrade, Serbia}
\affiliation{$^4$Brunel University London,  Uxbridge Middlesex, UB8
  3PH, UK}
\affiliation{$^5$Department of Health Sciences Research, Center for
Individualized Medicine, Mayo Clinic, Rochester, Minnesota 55905, USA}
\affiliation{$^6$Complexity Science Hub Vienna, Josephstadter Strasse
    39, Vienna, Austria}
\vspace*{3mm}
\date{\today}

\begin{abstract}
\noindent 
Motivated by the relevance of higher-order interactions in quantum physics and materials
science at the nanoscale, recently a model has been introduced for new classes of networks that grow by the geometrically constrained aggregation of
simplexes (triangles, tetrahedra and higher-order cliques). 
Their key features are hyperbolic geometry and hierarchical architecture with simplicial
complexes, which can be described by the algebraic topology of graphs. Based on the model of chemically tunable
self-assembly of simplexes [\v{S}uvakov et al., Sci.Rep 8, 1987
(2018)], here we study the impact of defect simplexes on the course of
the process and their organisation in the grown nanonetworks
for varied chemical affinity parameter and the size of building simplexes.
Furthermore, we demonstrate how the presence of patterned defect
bonds can be utilised to alter the structure of the assembly after the growth process is completed. In this regard, we consider the
structure left by the removal of defect bonds and quantify the changes
in the structure of simplicial complexes as well as in the underlying
topological graph, representing 1-skeleton of the simplicial complex.  By
introducing new types of nanonetworks, these results open a promising
application of the network science for the design of complex
materials. They also provide a deeper understanding of the mechanisms
underlying the higher-order connectivity in many complex systems.
\\[3pt]
\end{abstract}

\maketitle
\section{Introduction}
In recent years, the application of graph theory to analyse complex
patterns in  empirical data has revolutionised research in various
interdisciplinary sciences. Some well-known examples include emotion-driven online dynamics with co-evolving networks of users and posts
\cite{we-Entropy2013} and mapping brain imaging data (see recent related  \cite{we-BrainSciRep2019} and references there).  However, the use of graphs for  mapping certain problems in
physics and materials science is still in its infancy
\cite{we_nanonetworks2013,AT_Jap_book2015,Geometries_QuantumPRE2015,Spectra-rapisarda2019,we_Qnets2016,AT_chemistry}.
In this case, more profound knowledge of the physics and chemistry of the problem helps to appropriately identify the nodes and edges of the structure, that often refers to the phase space of the system rather than a real-space structure.  

In materials science, complex structures made of nano-scale objects
often correlate  with an increased functionality
\cite{SA_NANOfocus2012,ProgMatt_review}. 
Processes of self-assembly are widely used  to
grow such systems, where the addition of each object to the growing structure obeys certain rules and locally minimises the energy \cite{SA_allscales_Science2002,minEclusters_PCCP2016,SA_colloids_scales2016}.
Therefore, the use of mathematical concepts \cite{AT_Materials_Jap2016} and graph representations
\cite{we_nanonetworks2013} are highly desirable for both the design and characterisation of the nanostructured assemblies. In this
context, real-space networks are visualised, for example, with the
nanoparticles as nodes and edges representing a kind of chemical binding \cite{we_JSTAT2009} or another association between them that is relevant for the problem in question. For example, the network representations of the conducting nanoparticle films have been studied in
\cite{we_NL2007,we_Oxford,we_topReview2010}. The use of graph theory
has enabled the description and differentiation of the structures that promote enhanced conduction via single-electron tunnellings between nanoparticles spaced within the quantum tunnellings radius in the direction of the current. 

Cooperative self-assembly \cite{SA_collectiveNR_SoftMatt2013,SA_hierarchical_SoftMatt2014,SA_Flowers}, where the pre-formatted group of
particles attach to a growing structure   represents a higher level
of the self-assembly processes, and opens an avenue towards new types
of \textit{materials inspired by mathematics} 
\cite{AT_Materials_Jap2016,we-SciRep2018}. Colloids with ``valence'' and directional bonding are a physical reality
\cite{FunctColloids_Nat2012}, particles with $n\in [1,7]$ active patches
were created by 2-stage swelling of the minimal moment clusters, and
subsequent DNA functionalization, resulting in different forms as
spheres, dumbbells, triangles, tetrahedra and higher-order structures.
In these processes, the geometrical-compatibility constraints of the binding forms with the growing structure play an important role, apart from the chemical affinity between the structure and the binding nanoparticles. Recently, we developed a model
with  the appropriate self-assembly rules \cite{we-SciRep2018,we-ClNets-applet}, where the
building blocks are suitably described by simplexes, i.e., edges, triangles, tetrahedrons, and cliques of higher orders.  
 A prominent feature of these structures is a hierarchical architecture of simplicial complexes, which is accessible to the methods of the algebraic topology of graphs \cite{jj-book,Q-analysis-book1982}, as well as emergent hyperbolicity in the graph's metric space \cite{HB-spaces-decomp2004,HB-Bermudo2016,HB-BermudoHBviasmallerGraph2013,Hyperbolicity_cliqueDecomposition2017}.

In this work, based on the model in   \cite{we-SciRep2018}, we extend the study of the cooperative self-assembly by considering the presence of defect simplexes and describe the impact of defect bonds on the assembled nanonetworks.  Specifically, we show how the presence of simplexes with a defect edge can alter the course of the process leading to the structure with not-random patterns of defect bonds and changed topological features of the assembly, depending on the size of the binding simplexes and the chemical affinity parameters. We further demonstrate how the patterned defect bonds can be utilised to alter the structure in the already grown assembly. Changes in the topological properties of the
assemblies are quantified using the algebraic topology analysis of simplicial complexes,  as well as by determining the hyperbolicity parameter of the underlying graph.

\begin{figure*}[!]
\begin{tabular}{cc} 
\resizebox{34pc}{!}{\includegraphics{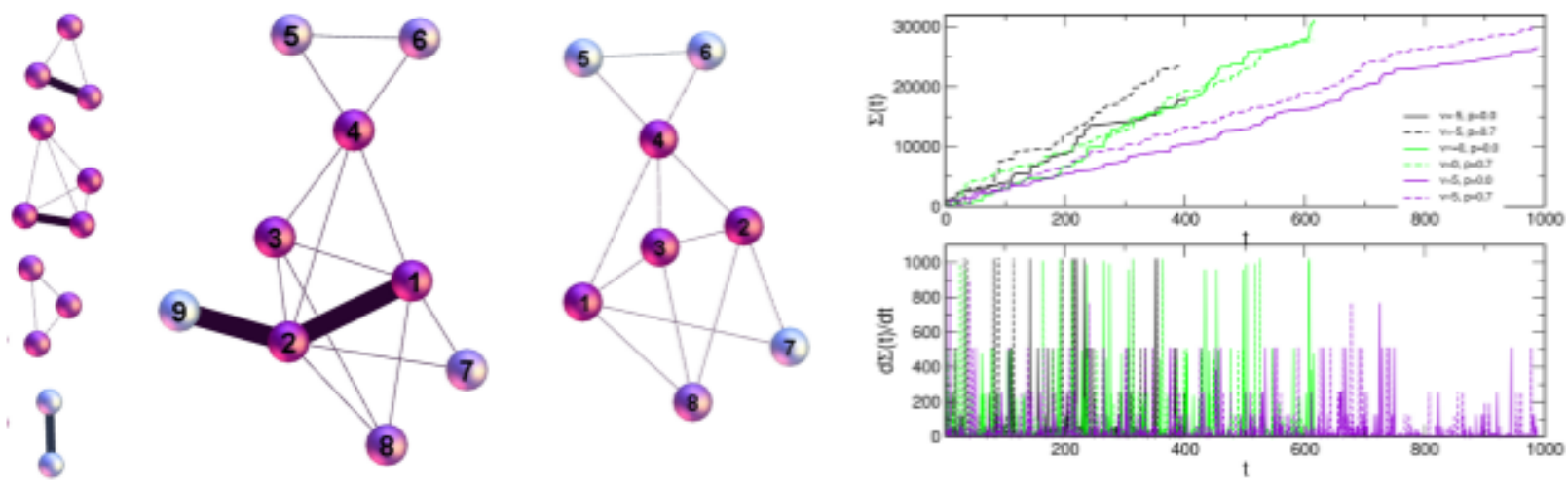}}\\
\end{tabular}
\caption{Left: Illustration of the aggregation of simplexes with defect
  edges, indicated by thick lines. For example, assuming that the
  structure shown in the middle emerged from  the type of simplexes shown on the
  left, we observe the following sequence of events: Starting
  from a defect tetrahedron with the vertices 1-2-3-4, a pure triangle
  4-5-6 is attached sharing the node 4, then a defect triangle 1-2-7
  is attached sharing the defect edge 1-2, following by the attachment
  of another defect tetrahedron 1-2-3-8 by sharing its defect triangle
  face with the previous structure; eventually, a defect dumbbell 9-2
  is attached sharing its defect face, the node 2. The connected structure that remains after removal of
  two defect edges is shown on the right of the figure.
Right: For the distribution of the attaching simplexes
  in the range $n\in [2,10]$, the evolution of the
  number of simplexes and faces 
in the graph $\Sigma (t)$  until the number of nodes reaches
1000. Three values of the affinity parameter $\nu$ are considered,
shown in the legend, combined  with
the probability $p=0.7$ of a
  defect bond in simplexes and simplexes with all pure bonds
  $p=0.0$. The lower panel shows the change
  in the number of simplexes in time.}
\label{fig-demo-Sigma} 
\end{figure*}
\begin{figure*}[!]
\begin{tabular}{cc} 
\resizebox{32pc}{!}{\includegraphics{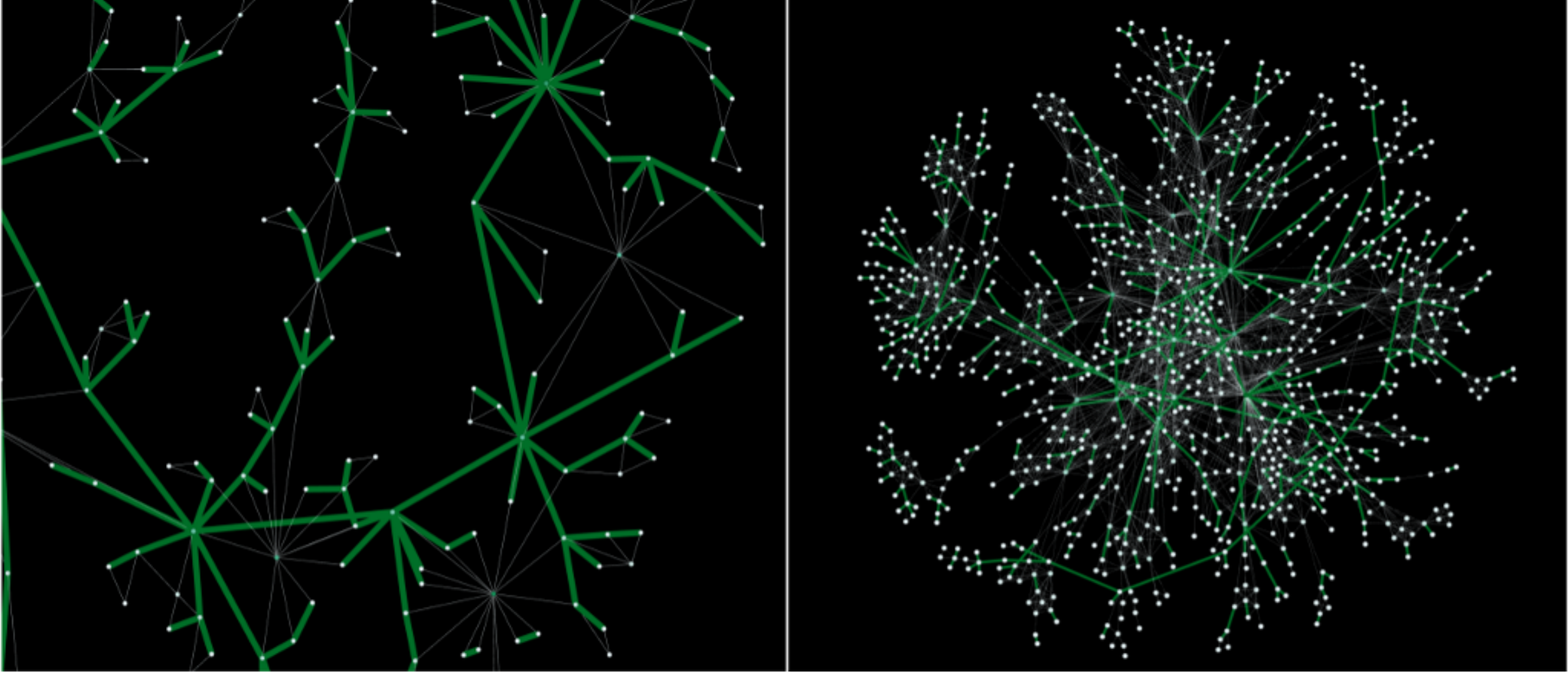}}\\
\end{tabular}
\caption{ Close-up of the structure of assembly of triangles  (left)
  and the assembly of the distributed clique sizes
  $n\in[2,10]$ according to $\sim n^{-2}$ (right) for strictly
  geometric aggregation ($\nu=0$) and  70\% defect
  simplexes. Defect edges are shown as thick (green) lines.}
\label{fig:netni0p07}
\end{figure*}

\section{The Model and growth of nanonetworks with defects\label{sec:model}}
Following the original model  \cite{we-SciRep2018}, we consider the assembly of simplexes that are full graphs (cliques) of $n$ vertices.
Starting from an initial simplex the new simplex  of  the size 
driven from the probability $p_n\sim n^{-\alpha}$ (we fix $\alpha =2$
if not specified) is attached to
the growing structure by sharing one of its faces. Notice that a
simplex of $n$ vertices possesses faces as sub-simplexes of all orders $q=0,1,2,3\cdots q_{max}-1$ from the vertex to the largest
sub-simplex,  where $q_{max}=n-1$ indicates the order of the
simplex. Determining the face to be shared depends first on the number
of geometrically compatible sites in the current structure
(geometrical compatibility constraint).  Moreover, when a face of the
order $q$ is shared with an already existing simplex in the growing
structure, the remaining $n_a=q_{max}-q$ vertices will be added to the
system as a formated group of nanoparticles. The affinity of the system towards the addition of a group is described by the chemical affinity parameter $\nu$
\cite{we-SciRep2018,we_SAloops}, which modifies the probability
defined based on the geometry factor alone, see Eq.\
(\ref{eq-pattach}).  Specifically, for a large negative $\nu$, the
system likes the addition of particles, which results in sharing a
minimal face, that is a single node. In this limit, the cliques are
effectively repelling each other. Whereas, in the opposite limit, with
a large  $\nu >0$, a single node is preferably added; thus, the added
clique shares its largest face with a previous compatible structure \cite{we-SciRep2018}.  

\begin{figure*}[!htb]
\begin{tabular}{cc} 
\resizebox{32pc}{!}{\includegraphics{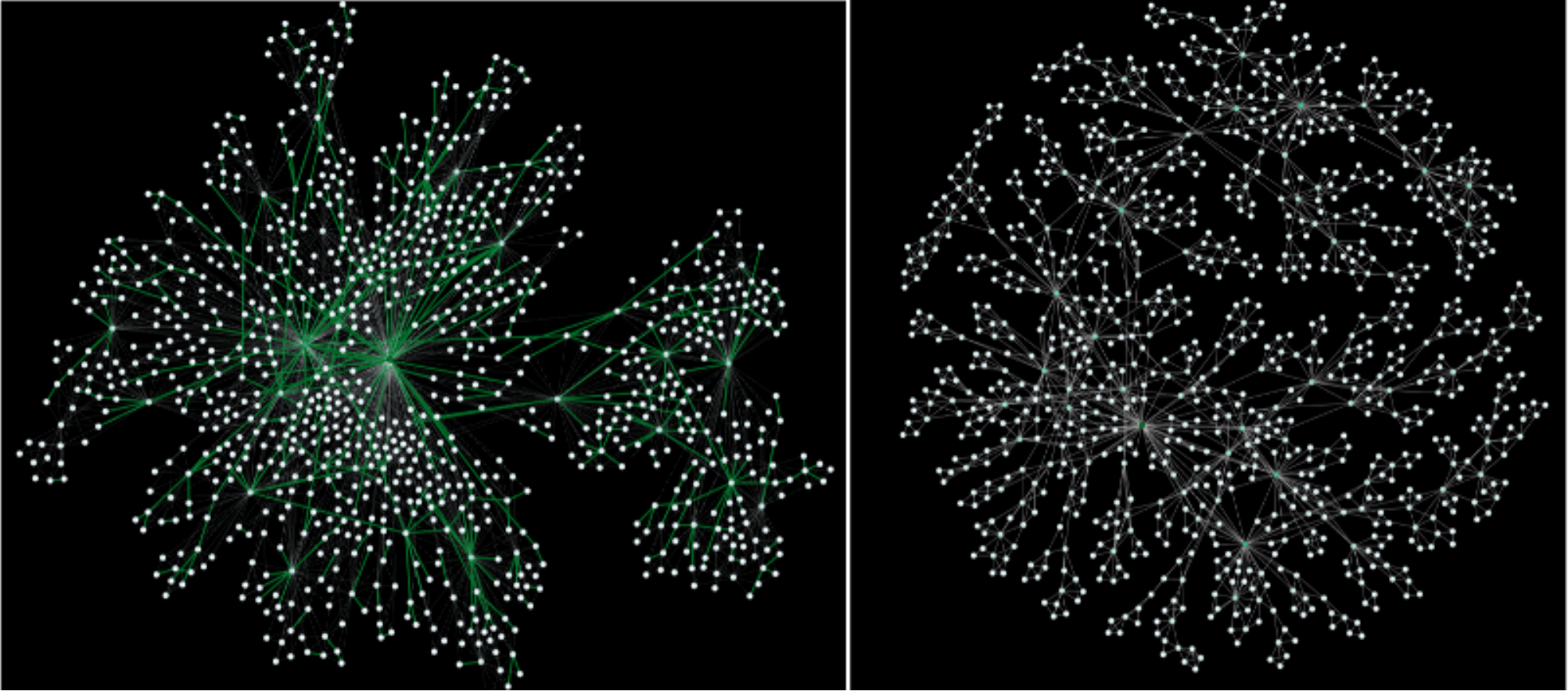}}\\
\end{tabular}
\caption{An example of the nanonetwork with defect bonds (thick green edges)
 assembled at affinity $\nu =$ +5 from tetrahedra with the
 probability of a defect bond $p=0.5$ (left),   and the structure
 obtained after removing the defect bonds (right).} 
\label{fig:tetra_p05ni5_rmw2}
\end{figure*}

Here, we allow that a simplex to be added to the structure can, with a
finite probability $p$, have a bond that differs from the other bonds,
i.e., a defect; the presence of such defect bonds affects the geometrical compatibility factors, as explained below.
More precisely, we have
\begin{equation}
p(q_{max},q;p,t)= \frac{c_q(p,t)e^{-\nu (q_{max}-q)}}{\sum _{q=0}^{q_{max}-1}c_q(p,t)e^{-\nu (q_{max}-q)}} \ 
\label{eq-pattach}
\end{equation}
which defines the normalised probability that a clique of the order
$q_{max}$ will attach along its face of the order $q$, subject of the
presence of a defect edge. Specifically, a defect edge of the arriving
clique can nest along another defect edge on the network, else the
adjacent nodes are shared, while the straight (pure) edges align along
the straight edges in the geometrically compatible shapes. Therefore,
at the evolution time $t$, the number of the geometrically similar
docking sites $c_q(p,t)$ of the searched order $q$ also depends on the
potential defect edge to match the defect in the face of the newly
added clique.  In this way, the presence of defects that are already
built in the structure affects future binding events. Fig.\
\ref{fig-demo-Sigma} illustrates the impact of the presence of defect
edges in the process of self-assembly.  Note that the number of
geometrically compatible sites for docking along with the faces with
pure bonds changes even if a small number of  defect bonds are
present. For example, to add a new pure triangle along its largest
face, i.e., an edge, to the structure in the middle, apart from the
tree bonds in the top triangle, we can have only two more candidates,
the bonds 3-4 and 3-8.  In contrast, 15 bonds are available in the
case when the same structure consists of only pure bonds. Some
examples of grown assemblies with varied parameters are shown in Fig.\ \ref{fig:netni0p07} and in Fig.\ \ref{fig:tetra_p05ni5_rmw2}.

The impact of defect edges for a given $p$ also depends on the size and the dispersion of the attaching cliques, and the binding affinity $\nu$. When the affinity among cliques is significant,
such that they intend to share their maximal faces, the aggregation of
defect edges is more effective (see Fig.\ \ref{fig:tetra_p05ni5_rmw2}),
leading to stronger aggregation of defect bonds and constraints to the
remaining structure.  Notably, defect bonds make a particular pattern.
These effects are especially  pronounced in the case of small
cliques, where a defect edge provides a more severe restriction on the binding of the remaining faces. In the case of purely geometrical aggregation, $\nu=0$, the defect edges at sufficiently large concentration form tree-like structures and  "highways'' through the graph.
Consequently, the grown network with defect simplexes is different from the case when the simplexes with equal edges were used (i.e., $p=0$), see Table\ \ref{tab-properties-nu-p}. In the following, we employ  Q-analysis \cite{cliquecomplexes,Q-analysis-book1982,SM-book2017} to quantitatively describe the organisation of simplicial complexes in various aggregates grown in the presence of defect bonds, in comparison with the case $p=0$. We also analyse the changes in the
structures caused by the removal of the defect bonds.

\section{Q-analysis \& Impact of defect bonds on the architecture of
  nanonetworks\label{sec:structure}}
In this context, a simplex of the order $q_{max}=n-1$ is a full graph of $n$ vertices.  In a simplicial complex, two simplexes are $q$-connected if they share a face of the order $q$, i.e., they have at least $q+1$ shared nodes. The dimension of the considered simplicial complex equals the dimension of
the largest clique $q_{max}+1$ belonging to that complex. 
To describe the structure of simplicial complexes 
at different topology levels $q=0,1,2 \cdots q_{max}$, Q-analysis uses
notation from the algebraic topology of graphs 
\cite{atkin1976,Q-analysis-JJ,Q-analysis-book1982}. 
Specifically: \begin{itemize}
\item the first structure vector (FSV) components $\{Q_q\}$
denote  the number of $q$-connected components;
\item  the second structure vector  (SSV) components $\{n_q\}$
  corresponds to the number of simplexes of the order greater than or equal to $q$;
\item  then the third structure vector  (TSV) component $q$ is determined as $\hat{Q_q}\equiv
1-Q_q/n_q$ measuring the degree of connectivity at the
topology level $q$ among the simplexes of the order higher than $q$.
\end{itemize} 
Using the Bron-Kerbosch algorithm \cite{BK} we construct the incidence
matrix $\Lambda (G)$ of the graph $G$, starting from its adjacency matrix.
The incidence matrix contains complete information about all present simplexes as well as the vertices that belong to each simplex. Thus the components of these structure vectors can be
determined from the corresponding incidence matrix  $\Lambda
(G)$. Further characterisation of the architecture of simplicial complexes is provided by the quantity $f_q$, which is defined \cite{we_PRE2015} as the number of simplexes and faces at the topology level $q$.  

Furthermore, we compute these topology features for the assembly that
remains after the removal of defect bonds. Notice that the removal of
the defect bonds changes the structure of the assembly by breaking the
simplexes in which such bonds were built-in. For example, cf.\  Fig.\ref{fig-demo-Sigma}, a defect tetrahedron breaks into two triangles that are  attached along the common edge when the defect bond is removed. The effects of the defect bond removal correlate with the size of the original cliques. More precisely, the following rule applies:
\begin{itemize}
\item a clique of the order $q_{max}$ with a defect bond breaks into two cliques of the order $q_{max}-1$; 
\item these new cliques are attached along their largest face, that is they share a face of the order $q_{max}-2$. 
\end{itemize}

\begin{figure*}[!]
\begin{tabular}{ccc} 
\resizebox{24pc}{!}{\includegraphics{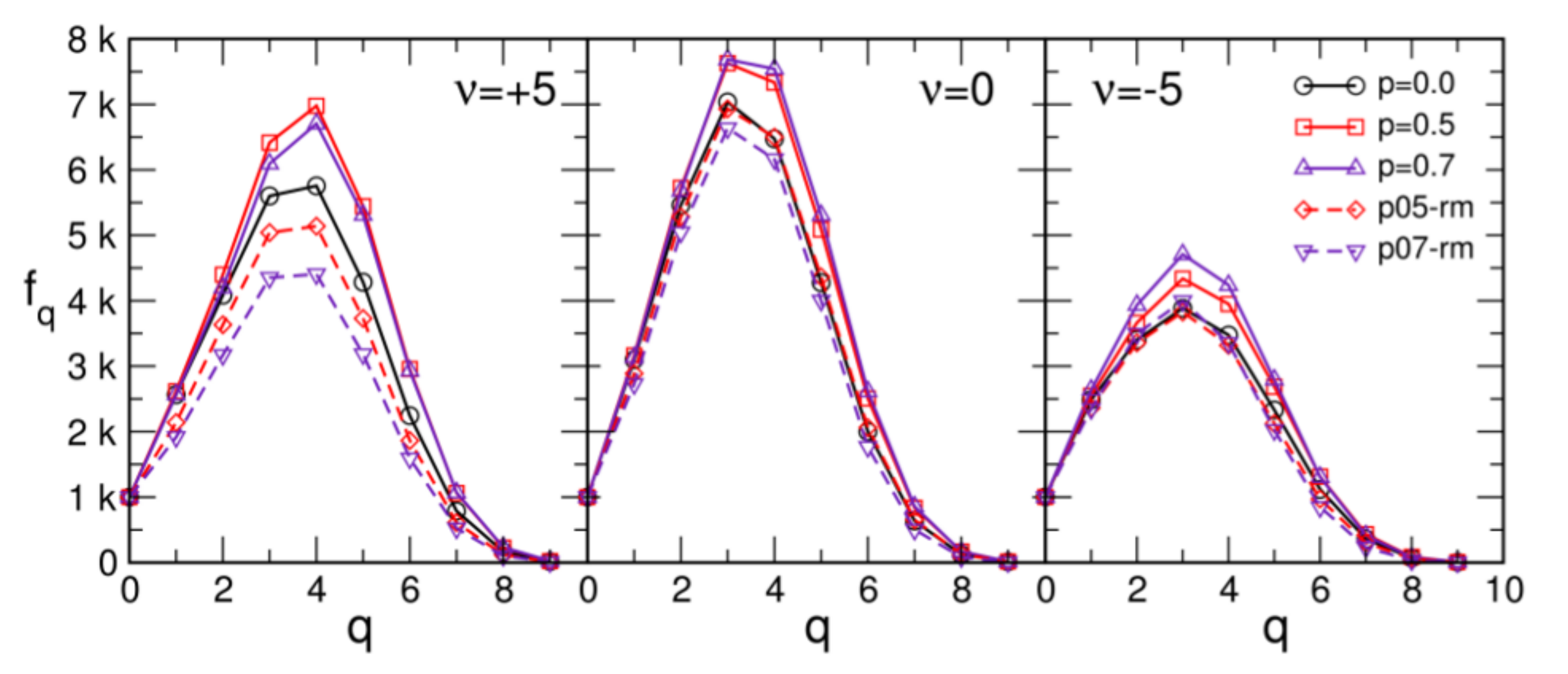}}\\
\end{tabular}
\caption{$f_q$ vs $q$ for the pure ($p=0$) and defect network
  ($p$=0.7, 0.5)
  and the network obtained by  removal of defect bonds ($p$=0.7- rmw2,
  p05-rmw2)
  for three values of $\nu=$5,0,-5, left to right.}
\label{fig:fq3x3}
\end{figure*}

\begin{figure*}[!htb]
\begin{tabular}{ccc}
\resizebox{24pc}{!}{\includegraphics{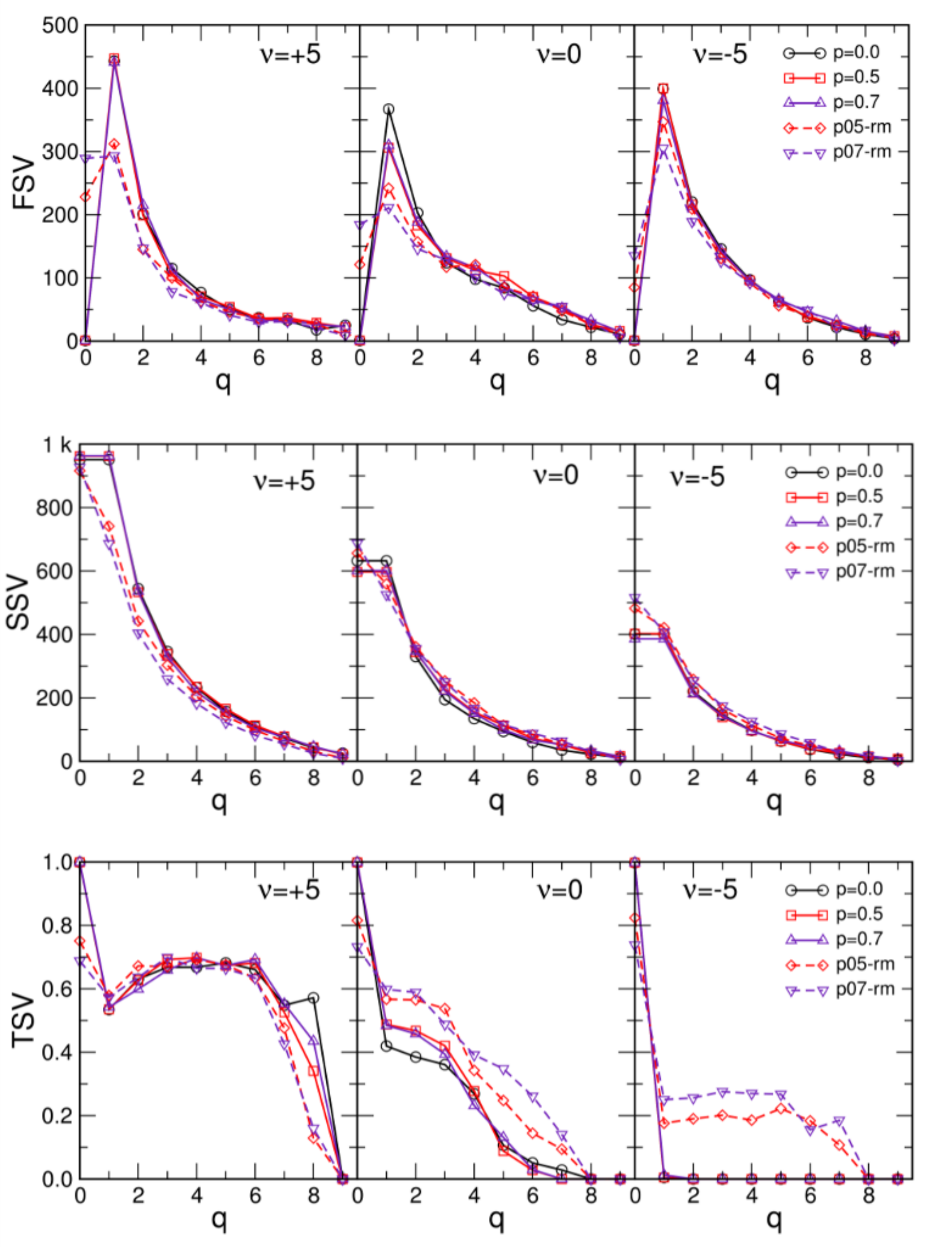}}\\
\end{tabular}
\caption{Components of the first (FSV), second (SSV)  and third (TSV) structure
  vector (top-to-bottom row) against the topology level $q$ for the pure ($p$=0.0) and
  defect network ($p$=0.5 and 0.7)
  and the network obtained by  removal of defect bonds (p05-rmw2, p0.7- rmw2)
  for three values of $\nu=$5,0,-5, indicated on each panel.}
\label{fig:SVs3x3}
\end{figure*}

 Fig.\ \ref{fig:fq3x3} and Fig.\ref{fig:SVs3x3} show these structural
 properties of a  few representative aggregates of simplexes obtained
 with/without defect bonds. Notably, the aggregation of defect faces
 causes those made of pure bonds to spread, which results in higher
 values of $f_q$ for a finite concentration of defect bonds $p$, as
 compared with $p=0$. This principle applies,  although the values are
 different, for all aggregates at different values of parameter
 $\nu$. One should also notice that the highest point  of
$f_q$ correlates with the number of simplexes that need to be added to
complete a given number of vertices, here $N\geq 1000$, which is considerably
different for different $\nu$, cf.\ Fig.\ \ref{fig-demo-Sigma}, right.
With the removal of the defect bonds, generally, we have a smaller number of simplexes and faces, resulting in the proportional decrease of $f_q$ at all $q$ levels, in comparison with the original
structure with defects. The effects are more pronounced
in the dense structure of cliques, corresponding to $\nu >0$, cf.\
Fig.\ \ref{fig:fq3x3}, than in the structures with sparsely connected cliques (for $\nu<0$ and partly $\nu=0$).

In the structure vectors, shown in Fig.\ \ref{fig:SVs3x3}, notice that the aggregates grown with or without defect represent one connected component, then the FSV component $Q_0=1$. The peak in the FSV at $q=1$ and the decay at larger $q<10$ reflects the actual distribution of the size $n\in[2,10]$ of the attaching cliques; the distribution  $f_n\sim n^{-2}$ favours dumbbells as compared with higher-order cliques up to the 10-cliques. 
The total number of cliques, which is given by the 0th components of the SSV, given in the middle row in Fig.\ \ref{fig:SVs3x3}, is
significantly more prominent in compact structures ($\nu=+5$) than in
the sparse assembly of cliques grown at $\nu=0$ and $\nu=-5$. This
observation is compatible with the growth process depicted in Fig.\
\ref{fig-demo-Sigma}. 
In the presence of defect cliques, the number of $q$-connected components, as well as the total number of cliques from the level $q$ upwards, differ from the case of a structure with pure simplexes, in particular in the range of large and intermediate $q$-values. 
The components of TSV indicate that in the compact regions of the
graph (i.e., in the case $\nu=+5$ and partly in $\nu=0$), the large
cliques with defect bonds are more weakly connected than in the
pure-simplexes structure. 
However, in the range $q\in[2,6]$ for
$\nu=+5$ and $q\in [1,4]$ for the case $\nu=0$, the connectivity
exceeds the curve of TSV for the structure without defects.
Meanwhile, in the case $\nu=-5$, the cliques of all sizes share a
single node, therefore they appear to be disconnected already at the level $q=1$, cf.\ TSV in the lower right panel in Fig.\
\ref{fig:SVs3x3}.

With the removal of defect bonds, the number of large cliques gradually decreases, while the number of intermediate and small cliques results as a balance between breaking the initially present defect cliques of that order and the appearance of new ones from the broken defect cliques of 
one order higher.
Consequently, the FVS changes such that $Q_0$ increases because of broken bonds, some separate graph parts can occur. 
The changes are most dramatic in the case of $\nu<0$. Following a broken bond in a clique of order nine, we have two cliques of the order eight that are sharing a clique of the order seven,
and so on, as explained above. Consequently, non-trivial connectivity
appears among these newly generated cliques at all levels $q\in[1,8]$, 
as shown in the lower right panel
of Fig.\ \ref{fig:SVs3x3}, even though the originally
built-in cliques repelled each other such to share a single vertex. A
similar effect occurs in the sparse areas of the structure grown in
the absence of chemical factors ($\nu=0$). 
The effects are proportional to the probability of a defect bond $p$,
which decides the actual  number of defect bonds in the grown structure, depending on $\nu$ (see Table\ \ref{tab-properties-nu-p}).
In the following, we analyse how the changed architecture of
simplicial complexes due to breaking defect bonds affects the
hyperbolicity and other features of the topological graph.

\section{Changes of hyperbolicity induced by the removal of defect
  bonds\label{sec:structure}}
As mentioned in the Introduction, the assembly of cliques possesses a \textit{negative curvature in the graph's metric space}, which is endowed with the shortest-path distance.  Hence, the generalised Gromov's 4-point hyperbolicity criterion can be applied to characterise it.
Specifically, the graph $G$ is hyperbolic \textit{iff} there is a
constant $\delta(G)$ such that for any four vertices $\{ A, B, C, D\}$
of the graph, the relationships between the sums of distances between distinct pairs of these nodes 
$d(A, B) + d(C, D) \le d(A, C) + d(B, D) \le d(A, D) + d(B, C)$
implies that 
\begin{equation}
\delta(A,B,C,D) = \frac{{\cal{L}}-{\cal{M}}}{2} \le \delta.
\label{eq-hyperbolicity}
\end{equation}
Here $d(U,V)$ indicates the shortest-path distance and we denoted the
largest  ${\cal{L}} = d(A, D) + d(B, C)$ and the middle value  ${\cal{M}} = d(A, C) + d(B, D)$.  
Observing that the upper bound of the expression in (\ref{eq-hyperbolicity}) is $({\cal{L}}-{\cal{M}})/2\le
d_{min}$, where  $d_{min}=min\lbrace d(A,B),d(C,D)\rbrace$ enables us to
determine the hyperbolicity parameter $\delta(G)$ by plotting
$({\cal{L}}-{\cal{M}})/2$ against $d_{min}$ and  investigating the worst case growth of the dependence. 
For each graph, using its adjacency matrix, we first compute the
matrix of distances between all pairs of nodes. Then, by sampling a
large number of sets of nodes for the 4-point condition
(\ref{eq-hyperbolicity}) we determine and plot the  largest $\delta $
against  the corresponding distance $d_{min}$; the
maximum observed value of all $\delta_{max}$  determines the graph's
hyperbolicity parameter $\delta (G)$.

\begin{figure*}[!hbt]
\begin{tabular}{cc} 
\resizebox{24pc}{!}{\includegraphics{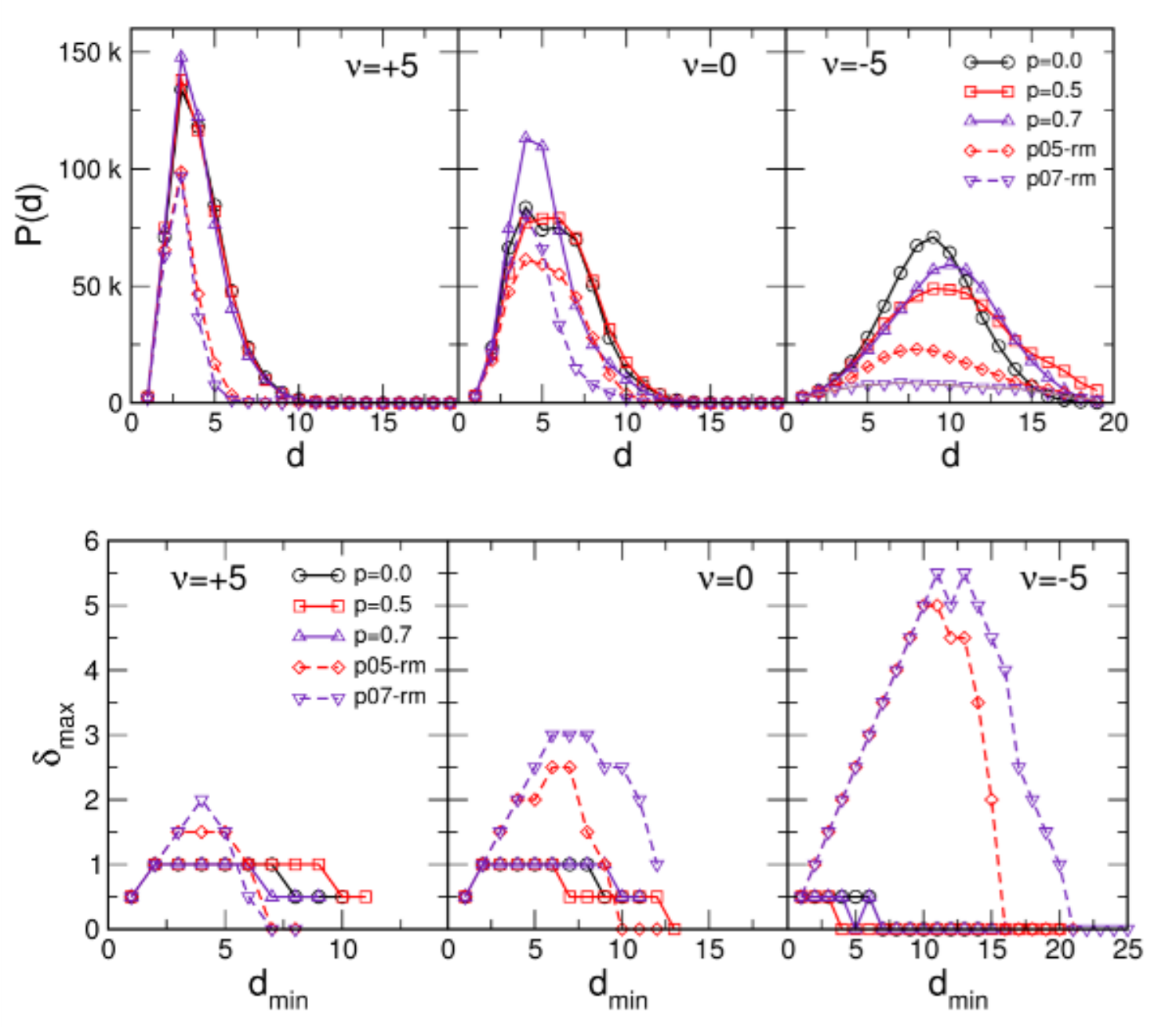}}\\
\end{tabular}
\caption{The distance distribution $P(d)$ vs $d$, and $\delta_{max}$ vs $d_{min}$ for the pure ($p=0$) and defect network ($p$=0.7,0.5)
  and the network obtained by  removal of defect bonds (p05- rmw2,p07-rmw2)
  for three values of $\nu=$5,0,-5, left to right.}
\label{fig:HB-Pd3x3}
\end{figure*}
For the graphs with a small hyperbolicity
parameter, it is known \cite{HB-Bermudo2016,HB-distortionDelta2014,Hyperbolicity_cliqueDecomposition2017} that
the upper bound of the hyperbolicity parameter is related to a specific subjacent
structure; for example \cite{Hyperbolicity_chordality2011}, the presence of an isometric cycle  $C_n$ of
the length $n\geq 3$ would lead to $\delta (C_n)=\lfloor
  n/4\lfloor - \frac{1}{2}$, if $n\equiv 1 (mod 4)$, else $\delta (C_n)=\lfloor
  n/4\lfloor$.  Similarly,
since the cliques are ideally hyperbolic ($\delta_{clique}=0$), a
combination of cliques that are apart at a small distance $i$ causes the
increase of the hyperbolicity parameter by an integer
\cite{Hyperbolicity_cliqueDecomposition2017}, i.e.,
$\delta_{clique}+i$.

\begin{table*}
\caption{Graph measures of the assemblies of simplexes of the size $n\in
  [2,10]$ distributed as $\varpropto n^{-2}$ and the probability of
  defect bond $p$,
  for three representative values of the affinity parameter $\nu=\pm$
  5, and 0. The properties of graphs with removed pattern of defect edges, (0.7-db)
  and the graph when the same number of edges $c$ is removed at random
  (rand-c) are also shown. 
  The effective concentration of defect edges $c$, the average degree $\langle k\rangle$, path length  $\langle \ell \rangle$ and clustering coefficient $\langle
  Cc\rangle$, graph's modularity $mod$, diameter $D$, spectral
  dimension $d_s$, all computed for the graph size $N=5000$ nodes. Additional
  properties computed for the graphs of the same parameters but
  $N=1000$ nodes are the hyperbolicity parameter
  $\delta_{max} $, and the topology level $q^*$ at which the
  connectivity (third structure vector TSV) between the simplexes drops to
  zero, and the connectivity at the level before $q^*-1$. }
\begin{tabular}{|c|c|cccccc|c|c|c|}
\hline\hline
$\nu$  &$p$&$c$& $<k>$&$<\ell>$&
                                 $<Cc>$&mod&D&$\delta_{max}$&$q^*$&$TSV(q^*-1)$\\
\hline
  & 0.0 & 0& 5.005& 4.475&0.601& 0.524&18& 1&  9&0.057\\
+5&0.7&0.271&5.115&4.197&0.602&0.556&19&1& 9&0.435\\
& 0.7-db& 0&5.0671&3.025&0.774&0.414&11&2.0&9&0.160\\
&rand-c&0&4.162&3.971&0.492&0.517&14&3.0&6&0.079\\
\hline
& 0.0 &0& 5.988& 6.209& 0.714&0.882&17&1&8&0.0285\\
0&0.7&0.149&5.933&6.256&0.721&0.872&17&1&7&0.0298\\
& 0.7-db&0&6.124&5.719&0.742&0.867&17&3.0&8&0.141\\
&rand-c&0&5.231&7.027&0.610&0.883&23&3.0&7&0.024\\
\hline
& 0.0 &0& 5.075& 13.213&0.813&0.972&32&1&2&0.005\\
-5&0.7&0.109& 5.270&11.788&0.825&0.966&27&1&2&0.0129\\
& 0.7-db&0&5.223&10.417&0.783&0.975&31&5.5&8&0.185\\
&rand-c&0&4.625&14.751&0.730&0.973&31&4.5&7&0.218\\
\hline
\hline
\end{tabular}
\label{tab-properties-nu-p}
\end{table*}
The network growth in our algorithm  by attaching a new clique such
that it shares a face with another previously present clique in the system, immediately implies that their hyperbolicity parameter cannot exceed unity. That is, these are 1-hyperbolic graphs
\cite{we-SciRep2018}, as also confirmed by a direct computation, see Fig.\ \ref{fig:HB-Pd3x3}. The same conclusion also applies to the structure grown with the defect cliques, as long as the cliques are complete. However, by removing the defect bonds, 
the cliques that contained them break into smaller cliques that appear
to be differently attached  to the rest of the graph. Moreover, holes of
different sizes and dimensions can occur. Consequently, we have
an increase of the hyperbolicity parameter of the whole graph, due to the
presence of longer cycles and increased distances between the newly
observed cliques. 

Fig.\ \ref{fig:HB-Pd3x3} shows the distribution of distances in the
case of pure simplexes and in the presence of defect cliques, and how it changes by the removal of defect bonds for varied parameters $p$ and $\nu$.  Notice that the distribution of distances between pairs of nodes
changes due to the presence of defect bonds.  
In contrast to the dense graphs (for $\nu=+5$), where the most
probable distance remains three, in the sparse graphs (at $\nu=0$ and
especially at $\nu=-5$) the most probable distances are larger  than in the case without defects. Similarly, these graphs experience the most dramatic changes in the distance distributions when the defect bonds are removed.  The diameter of the graph (referring to the largest connected component) also 
changes, cf.\ table \ref{tab-properties-nu-p}.

The lower panels show the  hyperbolicity parameter
$\delta_{max}$ for the corresponding graphs. As expected,
$\delta_{max}=1$ for all graphs grown by the attachment of cliques
rule with and without defect bonds, for all $\nu$ values. However,
when the defect bonds are removed, the changed organisation of
simplexes, as described above, leads to the appearance of holes and
long cycles, which results in the increased values of
$\delta_{max}$. The increase strongly depends on $\nu$ at which the
graph with defect simplexes is grown. More precisely, in compact structures, the hyperbolicity parameter reaches the values 3/2 or 2, compatible with the subjacent structures with increased distances between the cliques.
On the other hand, a substantial increase of $\delta_{max}$
 in the sparse structures, reaching the values 3 for $\nu=0$ and 5.5
 for $\nu=-5$  and the considered concentration of the removed bonds,  can be related to the appearance of long cycles.

\section{Summary and discussion\label{sec:discussion}}
We have introduced a model for self-assembly of simplexes in the
presence of a defect bond. The model allows for the variation of the
probability of a defect bond in the attaching simplexes in conjunction
with their  size $n$ and the chemical affinity $\nu$, leading to a rich variety of resulting assemblies.
In this study, we consider a fixed probability $p$ and the sizes distributed according to $p_n\sim n^{-2}$  in
the range $n\in[2,10]$. We have shown how the presence of defect bonds can tune the structure of simplicial complexes as well as the underlying topological graph. The results of the quantitative analysis in  Fig.\
\ref{fig:fq3x3}, Fig.\ \ref{fig:SVs3x3} and Fig.\ \ref{fig:HB-Pd3x3},
show that the model provides the framework to grow a rich structure of
simplicial complexes with the possibility to control both the process
of the growth of the assembly as well as to change it by influencing
the defect edges after the growth is completed. In this study, we have demonstrated how the removal of defect
bonds leads to the altered structure of simplicial complexes;
moreover, the presence of holes and cycles in these transformed
assemblies are associated with  the increase of the graphs hyperbolicity parameter. 
Some standard  graph properties and their hyperbolicity as well as
measures indicating the connectivity between simplexes are listed in
the table \ref{tab-properties-nu-p} for the representative sets of parameters.

Remarkably, the defect bonds make non-random patterns---tree-like
structures in the graphs, even though no long-range forces are
present.  The apparent attraction among defects is primarily related to the geometrical constraints for the docking of simplexes; thus, it depends on the size of simplexes and the chemical affinity towards new vertices. Therefore, the removal of patterns of defect bonds has a profound effect on the structure of simplicial complexes, as discussed above.
These patterns also shape the graph's properties differently, as compared with the case when the same number of defect bonds are
randomly distributed, cf.\ table\ \ref{tab-properties-nu-p}.

In summary, we have introduced new classes of networks that evolve by
self-assembly of simplexes with different shapes and types of
bonds. The variations of the parameters governing the process of
self-assembly allow different types of structures to grow from
sparsely separated simplexes to compact structure with large
simplicial complexes, and the possibility to modify their organisation
by affecting a specific type of bonds. These approaches are suitable
for designing new classes of nano-structured assemblies and for their 
quantitative characterisation beyond the standard pairwise
interactions.   The presented study also offers a deeper understanding
 of the mechanisms beyond the higher-order connectivity that lead 
to the occurrence of simplicial complexes in many other complex systems, from human connectomes
\cite{we-BrainSciRep2019} to patterns representing the brain-to-brain
coordination \cite{we-Brain-PLOS2016} and online social dynamics
\cite{we_PhysA2015,we_Tags2016}, as well as a variety of  problems in
physics  \cite{we_PRE2015,we_Qnets2016,Geometry_colloidsPRE2015,Geometries_QuantumPRE2015,Spectra-rapisarda2019}.

\section*{Acknowledgments}

Authors acknowledge the financial support from the Slovenian
Research Agency under the program P1-0044
and  from the Ministry of Education, Science and Technological Development of
the Republic of Serbia, under the projects OI 171037,  III 41011 and OI 174014.

\section*{Author contributions statement}
B.T. and M.S. designed research, M.S. contributed program tools, M.A.,
B.T., M.S., G.J.R. performed simulations and analysed data, B.T. produced figures and wrote the manuscript, all authors reviewed the manuscript.

\end{document}